\documentclass[12pt]{amsart}
\usepackage{amsmath,amssymb,amsbsy,amsfonts,amsthm,latexsym,amsopn,amstext,
                            amsxtra,euscript,amscd}
\usepackage{url}
\usepackage{enumerate}
\usepackage{subfig,graphics,epsfig,color}
\usepackage{float}

\begin{document}

\newtheorem{theorem}{Theorem}
\newtheorem{lemma}[theorem]{Lemma}
\newtheorem{algol}{Algorithm}
\newtheorem{cor}[theorem]{Corollary}
\newtheorem{prop}[theorem]{Proposition}

\newcommand{\comm}[1]{\marginpar{%
\vskip-\baselineskip %raise the marginpar a bit
\raggedright\footnotesize
\itshape\hrule\smallskip#1\par\smallskip\hrule}}

%%%%%%%%%%%%%%%%%%%%%%%%%
% Alphabet calligraphic %
%%%%%%%%%%%%%%%%%%%%%%%%%
\def\cA{{\mathcal A}}
\def\cB{{\mathcal B}}
\def\cC{{\mathcal C}}
\def\cD{{\mathcal D}}
\def\cE{{\mathcal E}}
\def\cF{{\mathcal F}}
\def\cG{{\mathcal G}}
\def\cH{{\mathcal H}}
\def\cI{{\mathcal I}}
\def\cJ{{\mathcal J}}
\def\cK{{\mathcal K}}
\def\cL{{\mathcal L}}
\def\cM{{\mathcal M}}
\def\cN{{\mathcal N}}
\def\cO{{\mathcal O}}
\def\cP{{\mathcal P}}
\def\cQ{{\mathcal Q}}
\def\cR{{\mathcal R}}
\def\cS{{\mathcal S}}
\def\cT{{\mathcal T}}
\def\cU{{\mathcal U}}
\def\cV{{\mathcal V}}
\def\cW{{\mathcal W}}
\def\cX{{\mathcal X}}
\def\cY{{\mathcal Y}}
\def\cZ{{\mathcal Z}}

\def\C{\mathbb{C}}
\def\F{\mathbb{F}}
\def\K{\mathbb{K}}
\def\Z{\mathbb{Z}}
\def\R{\mathbb{R}}
\def\Q{\mathbb{Q}}
\def\N{\mathbb{N}}
\def\M{\textsf{M}}

\def\({\left(}
\def\){\right)}
\def\[{\left[}
\def\]{\right]}
\def\<{\langle}
\def\>{\rangle}

\def\e{e}

\def\eq{\e_q}
\def\fS{{\mathfrak S}}

\def\lcm{{\mathrm{lcm}}\,}

\def\fl#1{\left\lfloor#1\right\rfloor}
\def\rf#1{\left\lceil#1\right\rceil}
\def\mand{\qquad\mbox{and}\qquad}

\def\jt{\widetilde\jmath}
\def\ellmax{\ell_{\mathrm max}}
\def\llog{\log\log}
\def\ordl{\mathrm{ord}_\ell\,}
\def\cl{{\mathrm cl}}
\def\Gal{{\mathrm Gal}}
\def\rk{\mathrm{rk}}
\def\Xl{\cX_\ell}
\def\Xld{\cX_{\ell,d}}
\def\Xlf{\cX_{\ell,f}}
\def\Xlm{\cX_{\ell,m}}
\def \Kl{\cK_\ell}
\def\Mlrs {M_{\ell,r,s}}
\def\Jac{\operatorname{Jac}}
\def\Res{\operatorname{Res}}
%%%%%%%%%%%%%%%  Topmatter %%%%%%%%%%%%%%%%%%

\title[Demjanenko Matrix of Quotients of Fermat Curves]
{On the singularity %% and Rank 
of the Demjanenko matrix of quotients of 
Fermat curves}

\author{Francesc Fit{\'e}}  
\address{Fakult{\"a}t f{\"u}r Mathematik, Universit{\"a}t Bielefeld, P.O.Box 100131, D-33501 Bielefeld, Germany} 
\email{francesc.fite@gmail.com}

\author{Igor E.~Shparlinski} 
\address{Department of Pure Mathematics, University of 
New South Wales, Sydney, NSW 2052, Australia}
\email{igor.shparlinski@unsw.edu.au}

%\subjclass{11G07, 11T06, 11Y16}


\begin{abstract} Given a prime $\ell\geq 3$ and a positive integer $k \le  \ell-2$, one can define a matrix $D_{k,\ell}$, the so-called Demjanenko matrix, whose rank is equal to the dimension of the Hodge group of the Jacobian $\Jac(\cC_{k,\ell})$ of a certain quotient of the Fermat curve of exponent $\ell$. For a fixed $\ell$, the existence of $k$ for which $D_{k,\ell}$ is singular (equivalently, for which the rank of the Hodge group of $\Jac(\cC_{k,\ell})$ is not maximal) has been extensively studied in the literature. We provide an asymptotic formula for the number of such $k$ when $\ell$ tends to infinity.
\end{abstract}

\maketitle

\section{Introduction}
For a prime $\ell\geq 3$ and a positive integer $k \le  \ell-2$, define the set
$$
M_{k,\ell}:=\{j\in (\Z/\ell \Z)^*\,|\, \langle kj\rangle_\ell+\langle j\rangle_\ell < \ell \}\,,
$$
where, for $j\in (\Z /\ell \Z)^*$, we denote by $\langle j\rangle_\ell$ the unique  
integer representative of $j$ modulo $\ell$ in the range $1,\dots,\ell-1$. This is a set of cardinality $(\ell-1)/2$. Koblitz and Rohrlich~\cite{KR78} show that the subgroup
$$
W_{k,\ell}:=\{ w\in (\Z/\ell \Z)^*\,|\, wM_{k,\ell}=M_{k,\ell}\}
$$
of elements stabilizing $M_{k,\ell}$ has cardinality $3$ or $1$ depending on whether
the parameter~$k$ is a primitive cubic root of unity modulo $\ell$ or not. The {\it Demjanenko matrix\/}  
is then defined as 
$$
D_{k,\ell}:=\left(E_{k,\ell}(-c^{-1}a)-\frac{1}{2}\right)_{c,a\in M_{k,\ell}/W_{k,\ell}}
$$
where 
$$
E_{k,\ell}(a):=
\begin{cases} 0 & \text{if $a\in M_{k,\ell}$,}\\
1 & \text{if $a\not\in M_{k,\ell}$.}\\
\end{cases}
$$
Consider now the curve
$$
\cC_{k,\ell}: \quad V^\ell = U(U+1)^{\ell-k-1}\,.
$$
This is a curve of genus $(\ell-1)/2$ that may be obtained as a quotient of the {\it Fermat curve\/} $\mathcal F_\ell: Y^\ell = X^\ell + 1$ by a certain subgroup of automorphisms of $\cF_\ell$; we 
refer to~\cite{FGL} for details. The fact that the rank of $D_{k,\ell}$ coincides with the dimension of the Hodge group of the Jacobian of $\cC_{k,\ell}$ has been 
exploited in~\cite{FGL} to determine the distribution of Frobenius 
traces attached to $\cC_{k,\ell}$ when $D_{k,\ell}$ is non-singular. 

Let $\Kl$ denote the set of positive integers $k \le \ell-2$
for which $D_{k,\ell}$ is singular.
It is easy to see that for every $\ell\equiv 2 \pmod {3}$, the set~$\Kl$ is empty 
(see Lemma~\ref{lemma: charac}).
Lenstra has shown (see~\cite[p.~354]{Gre80}) that $\Kl$ is non-empty for 
every sufficiently large $\ell\equiv 7 \pmod {12}$. In this note, we give an asymptotic formula
for the cardinality of $\Kl$, which, in particular, shows that $\Kl$ is non-empty for 
an overwhelming majority primes $\ell \equiv 1 \pmod 3$. 

\begin{theorem}
\label{thm:sing} Let $\ell-1= 2^\alpha 3^\beta m$ for some 
integers $\alpha > 0$, $\beta \geq 0$ and~$m$ with $\gcd(m,6)=1$.
Then 
$$
\left|\# \Kl- \frac{1}{2 ^{2 \alpha+2}}\(1-\frac{1}{ 3 ^{2 \beta}}\)\ell\right| \leq 4 \beta ^2\sqrt \ell + \frac{33}{16}\,.
$$
\end{theorem}

The key result to prove Theorem~\ref{thm:sing} is the characterisation of the non-singularity of $D_{k,\ell}$ in terms of certain conditions on the multiplicative orders of $k$ and $k^2+k$ modulo  $\ell$ obtained in~\cite{FGL} (see Lemma~\ref{lemma: charac} below).

\begin{cor}
\label{cor:positive} Let $\ell-1= 2^\alpha 3^\beta m$ for some 
integers $\alpha, \beta > 0$ and~$m$ with $\gcd(m,6)=1$.
If $\ell  >441 \cdot 2^{4\alpha }\beta^4$ then $\# \Kl> 0$.
\end{cor}

The previous result can be verified by direct calculations for $\ell \le 7$, 
and for $\ell > 11$ it follows from the inequalities
$$
 \frac{1}{2 ^{2 \alpha+2}}\(1-\frac{1}{ 3 ^{2 \beta}}\)\ell \ge 
\frac{1}{ 2 ^{2 \alpha-1}9}\ell 
$$
and
$$
 4 \beta ^2\sqrt \ell + \frac{33}{16} <  \(4 + \frac{33}{16\sqrt{11}}\)  \beta ^2\sqrt \ell 
 \le  \(4 + \frac{2}{3}\)  \beta ^2\sqrt \ell.
$$

We can now obtain an explicit form of the observation of Lenstra.
\begin{cor}\label{corollary: 7mod12}
For every prime 
$\ell\equiv 7 \pmod {12}$ distinct from 
$7$ and $19$
we have $\#\Kl> 0$.
\end{cor}

This can be deduced from Corollary~\ref{cor:positive} in the following way. First note that $\alpha=1$. Observe that for 
$$
m\geq m_\beta:=\lceil 441\cdot 2^3\cdot 3^{-\beta}\beta^4 \rceil\,,
$$
we have that 
\begin{equation}
\label{eq:alpha1}
\ell=2\cdot 3^ \beta m
\end{equation}
satisfies the hypothesis of Corollary~\ref{cor:positive} and thus $\#K_\ell>0$. Since for $\beta \geq 18$, one has $m_\beta=1$, we can limit our search for primes $\ell\equiv 7 \pmod {12}$  with $\#\Kl= 0$ among the finite set of primes $\ell$ of the form~\eqref{eq:alpha1}
with 
$$
\beta\in\{1,\dots,17\} \mand m\leq m_\beta-1 \ \text{with}\ \gcd(m,6)=1.
$$
A computer search establishes that the only primes of this form are 
$$
7,19,163, 487, 1459, 39367, 86093443, 258280327\,.
$$
Among the above primes, we have $\#\Kl= 0$ only for $\ell=7,19$.

As we have mentioned, Lemma~\ref{lemma: charac} below immediately implies 
that if $\ell\not \equiv 1 \pmod 3$ then $\Kl = \emptyset$. This is consistent 
with the vanishing of the main term of Theorem~\ref{thm:sing} for $\beta = 0$.
We also use Corollary~\ref{cor:positive} to derive a bound 
on the density of primes~$\ell\equiv 1\pmod 3$ with $\Kl= \emptyset$.

\begin{theorem}
\label{thm:Kl empty} For  $x \ge 2$ 
there are at most $O(x^{3/4}(\log x)^3)$ primes $\ell\equiv 1\pmod 3$ with $\ell \le x$ and $\#\Kl= 0$.
\end{theorem}

We can not answer the question of whether there exist infinitely many primes  $\ell\equiv 1\pmod 3$ with $\#\Kl= 0$. However, we provide a reason to believe so. Indeed, standard heuristic arguments suggest 
that for any $\beta\geq 1$ and $m\geq 1$ with $\gcd(m,6)=1$ there are infinitely many primes of the form $\ell=  2^\alpha 3^\beta m +1$, with $\alpha>0$, and we now show that $\# \Kl=0$ for most of such primes.
To this aim, for  fixed integers $\beta \ge 0$ and $m\geq 1$ with  with $\gcd(m,6)=1$, 
we define $\cL_{\beta,m}$ to be the set of primes of the form $\ell=  2^\alpha 3^\beta m +1$, for some $\alpha > 0$, such that $\#\Kl>0$. Then we have the  following 
finiteness result. 

\begin{theorem}
\label{thm:vanish} 
For any fixed~$\beta\geq 0$ and $m\geq 1$ such that $(m,6)=1$, the set $\cL_{\beta,m}$ is finite. More precisely, if $\beta=0$, then $\#\cL_{\beta,m}=0$, and we have
$$
\# \cL_{\beta,m} =O\( 3^{2\beta}m^ 2/\beta \)
$$
for $\beta \ge 1$. 
\end{theorem}

%\begin{theorem}
%\label{thm:vanish} 
%We have $\# \cL_\beta=0$ for $\beta= 0,1,2,3$ and
%$$
%\# \cL_\beta =O\(3^{2\beta}/\beta\)
%$$
%for $\beta \ge 4$. 
%\end{theorem}

\section{Preparations}

Let $\ordl k$ denote the multiplicative order 
of $k$ modulo $\ell$. Also for a prime $p$ and 
an integer $m$ we denote by 
$\nu_p(m)$ the $p$-adic order of $m$, that is, 
the largest integer $\nu$ with $p^\nu \mid m$. 
Our main tool is the following characterisation 
of the elements of $\Kl$ given in~\cite{FGL}.
%% [Theorem 1.2, \cite{FGL}]

\begin{lemma}\label{lemma: charac}
\label{lem:cond} For a prime $\ell\geq 3$ and a positive integer $k \le \ell - 2$, we have $k\in \Kl$ if and only if the three following conditions hold:
\begin{itemize}
\item[(i)] $\ordl k \ne 3$;
\item[(ii)]$\nu_2 (\ordl k) = \nu_2(\ordl(-k^2 - k)) = 0$;
\item[(iii)] $\nu_3 (\ordl k) > \nu_3(\ordl(k^2 + k))$.
\end{itemize}
\end{lemma}

Now, let $\Xl$ denote the group of multiplicative characters 
modulo~$\ell$. 
Furthermore, let $\Xld$ denote the set of 
characters of order dividing $d$, that is, the set of 
characters $\chi\in \Xl$ such that 
$\chi^d=\chi_0$, where $\chi_0$ is the principal character, see~\cite{IwKow} for a background on  characters.
We also use  $\Xld^*$ to denote the set of non-principal characters 
of $\Xld$. Given $\chi\in \chi_\ell$, we extend it to $\F_\ell$ in the following way: if $\chi = \chi_0$ is principal, then set $\chi_0(0):=1$. Otherwise, set $\chi(0):=0$.

Since $\Xl$ is dual to the multiplicative group $\F_\ell^*$ 
of the finite field of $\ell$ elements, for any divisor $t\mid \ell-1$ and $u\in \F_\ell^*$, for $d =  (\ell-1)/t$ we have

\begin{equation}
\label{eq:CharDiv}
\frac{1}{d} \sum_{\chi\in \Xld} \chi(u) = \left\{\begin{array}{llll}
1, &\quad\text{if } u^t = 1,\\
0, &\quad\text{otherwise}. 
\end{array}
\right.
\end{equation}
%In particular, the sum on the left hand side of~\eqref{eq:CharDiv}
%is the characteristic function of the property $\ordl u \mid t$. 
%We now use this observation together with the inclusion exclusion 
%principle to  detect the property  $\ordl u = t$ via 
%multiplicative characters as follows.
%For $d =  (\ell-1)/t$ we have
%\begin{equation}
%\label{eq:CharOrd}
%\sum_{f\mid d} \frac{\mu(f)}{f} \sum_{\chi\in \Xlf} \chi(u) = \left\{\begin{array}{llll}
%1, &\quad\text{if }\ordl u = t,\\
%0, &\quad\text{otherwise},
%\end{array}
%\right.
%\end{equation}
%where, as usual, we use  $\mu(f)$ to denote the M{\"o}bius function of 
%an integer $f\ge 1$.

Finally, we recall the following special case of the Weil bound of character sums
(see~\cite[Theorem~11.23]{IwKow}).

\begin{lemma}
\label{lem:Weil}
For any polynomial $Q(X) \in \F_\ell[X]$ with  $N$ distinct zeros
in  the algebraic closure $\overline \F_\ell$  of $\F_\ell$ and
which is not a perfect $s$th power in $\overline \F_\ell[X]$ for an
integer  $s\ge 2$, 
 and  a nonprincipal character $\chi \in \Xl^*$ of order $s$, we have
$$
\left| \sum_{k \in \F_\ell}\chi\(Q(k)\)\right| \le (N-1) \ell^{1/2}.
$$
\end{lemma}

\section{Proof of Theorem~\ref{thm:sing}}

Since condition {\it (i)} of Lemma~\ref{lem:cond}
fails to hold for at most two integers $k \in [1, \ell-2]$ we have 
\begin{equation}
\label{eq:KK}
|\# \Kl - \# \Kl^*| \le 2
\end{equation}
where $\Kl^*$ is the set of 
integers $k \in [1, \ell-2]$ satisfying the conditions {\it (ii)}
and {\it (iii)} of Lemma~\ref{lem:cond}. 

Let $\zeta(u)$ be the characteristic function of the condition
$\nu_2 (\ordl u)  = 0$. This is equivalent to 
$$
\ordl u \mid 3^\beta m =(\ell-1)/2^\alpha.
$$ 
So, we see from~\eqref{eq:CharDiv} that
\begin{equation}
\label{eq:zeta}
\zeta(u) = \frac{1}{2^\alpha} \sum_{\chi\in  \cX_{\ell,2^\alpha}} \chi(u)
=\frac{1}{2^\alpha} + \frac{1}{2^\alpha} \sum_{\chi\in  \cX_{\ell,2^\alpha}^*} \chi(u).
\end{equation}

Furthermore, for a non-negative integer
$h$, let $\eta_h(u)$ be the characteristic function of the condition
$\nu_3 (\ordl u) = h$. 
This is equivalent to 
$$
\ordl u \mid 2^{\alpha} 3^h  m =\frac{\ell-1}{3^{\beta-h}}
\mand
\ordl u \nmid 2^{\alpha} 3^{h-1}   m  
=\frac{\ell-1}{3^{\beta-h+1}}. 
$$ 
So, we see from~\eqref{eq:CharDiv} that
\begin{equation}
\label{eq:eta}
\begin{split}
\eta_h(u) & = \frac{1}{3^{\beta-h}}
 \sum_{\chi\in \cX_{\ell,3^{\beta-h}}} \chi(u)
 - \frac{1}{3^{\beta-h+1}}
 \sum_{\chi\in \cX_{\ell,3^{\beta-h+1}}} \chi(u)\\
& = \frac{2+\vartheta_h}{3^{\beta-h+1}} \\
 & \qquad \qquad  +
 \frac{1}{3^{\beta-h}}
 \sum_{\chi\in \cX_{\ell,3^{\beta-h}}^*} \chi(u)
 - \frac{1}{3^{\beta-h+1}}
 \sum_{\chi\in \cX_{\ell,3^{\beta-h+1}}^*} \chi(u),
 \end{split}
\end{equation}
where for in the case of $h=0$ we define
$\cX_{\ell,3^{\beta+1}}^*=\emptyset$
and we also set $\vartheta_0=1$ and $\vartheta_h=0$ for $h \ge 1$.

Then we have 
\begin{equation}
\label{eq:K and M}
\# \Kl^* =  \sum_{k=1}^{\ell-2} \cB_{k,\ell}=\sum_{k\in \F_\ell} \cB_{k,\ell}-\cB_{0,\ell}-\cB_{-1,\ell}\,,
\end{equation}
where
$$
\cB_{k,\ell} :=  \zeta(k) \zeta(-k^2 - k)
\sum_{r=1}^\beta \sum_{s=0}^{r-1}\eta_r(k) \eta_s(k^2 + k).
$$
Examining the expressions~\eqref{eq:zeta} and~\eqref{eq:eta} we conclude that 
each product $ \zeta(k) \zeta(-k^2 - k) \eta_r(k) \eta_s(k^2 + k)$, after expanding,
contains the constant term 
$$
\frac{1}{2^{2 \alpha}}\cdot \frac{2+\vartheta_r}{3^{\beta-r-1}} \cdot \frac{2+\vartheta_s}{3^{\beta-s-1}} = \frac{1}{2^{2 \alpha}}\cdot \frac{2}{3^{\beta-r-1}} \cdot \frac{2+\vartheta_s}{3^{\beta-s-1}}\,,
$$ 
(provided that $r \ge 1$ in our settings)
which does not depend on $k$, and also several terms with 
products of the form
$$
\chi_1(k) \chi_2(-k^2 - k)\chi_3(k) \chi_4(k^2 + k)
$$ 
with some characters 
$\chi_1,\chi_2\in \cX_{\ell,2^\alpha}$, $
\chi_3,\chi_4\in \cX_{\ell,3^{\beta-h}} \cup \cX_{\ell,3^{\beta-h+1}}$
such that at least one of them is nonprincipal. Since multiplicative 
characters form a cyclic group (see~\cite{IwKow}), we see that for some 
character $\chi$ of order $\ell-1$ and integers $f$, $g$, $h$
with $0\leq f,g <\ell-1$, $f+g >0$, $h = 0,1$ we have 
$$
\chi_1(k) \chi_2(-k^2 - k)\chi_3(k) \chi_4(k^2 + k) = 
\chi\(k^f(k^2+k)^{g}(-1)^{h}\)\,.
$$

Consider the polynomial
$$
P(\overline X, \overline Y, \overline Z, \overline T, \overline U, \overline V):= P_1(\overline X)P_2(\overline Y)\sum_{r=1}^\beta\sum_{s=0}^{r-1} P_{3,r}(\overline Z,\overline T) 
P_{4,r,s}(\overline U,\overline V)\,,
$$
where

\begin{equation*}
\begin{split}
& P_1(\overline X)=: \frac{1}{2^\alpha}+\frac{1}{2^\alpha}\sum_{i=1}^{2^\alpha-1}X_i\,,\qquad
P_2(\overline Y):=\frac{1}{2^\alpha}+\frac{1}{2^\alpha}\sum_{i=1}^{2^\alpha-1}Y_i\,,\\
&P_{3,r}(\overline Z, \overline T):=\frac{2}{3^{\beta-r+1}}+
 \frac{1}{3^{\beta-r}}
 \sum_{i=1}^{3^{\beta-r}-1} Z_{i,r}
 - \frac{1}{3^{\beta-r+1}}
 \sum_{i=1}^{3^{\beta-r+1}-1}T_{i,r}\,,\\
&P_{4,r,s}(\overline U,\overline V):= \frac{2+\vartheta_s}{3^{\beta-s+1}}+
 \frac{1}{3^{\beta-s}}
 \sum_{i=1}^{3^{\beta-s}-1} U_{i,r,s}
 - \frac{1}{3^{\beta-s+1}}
 \sum_{i=1}^{3^{\beta-s+1}-1}V_{i,r,s}\,,
\end{split}
\end{equation*}
and where $\overline X, \overline Y, \overline Z, \overline T, \overline U, \overline V$ are vector indeterminates given by
\begin{itemize}
\item $\overline X:=(X_i)$ and $\overline Y:= (Y_i)$ for $1\leq i \leq 2^{\alpha}-1$;
\item $\overline Z:= (Z_{i,r})$ and $\overline T:=(T_{i,r})$ for $1 \leq i\leq 3^{\beta-r}-1$ and $1 \leq r\leq \beta$;
\item $\overline U:= (U_{i,r,s})$ and $\overline V:=(V_{i,r,s})$ for $1 \leq i\leq 3^{\beta-s}-1$, $0 \leq s\leq r-1$, and $1 \leq r \leq \beta$.
\end{itemize}
Let $a_0, a_1,\dots, a_N$ be the set of coefficients of the polynomial $P$ with $a_0$ denoting the constant term. One observes that
\begin{equation}\label{eq: sumBk}
\sum_{k\in \F_\ell} \cB_{k,\ell}=\ell a_0+ \sum_{i=1}^N a_i \sum_{k\in \F_\ell} 
\chi_i\(k^{f_i}(k^2+k)^{g_i}(-1)^{h_i}\)\,,
\end{equation}
where for every $i=1,\dots,\ell-1$ we have that $0\leq f_i, g_i<\ell-1$, $f_i + g_i>0$ are integers,
$h_i =0,1$, and $\chi_i$ are characters of order $\ell-1$.  
Note that, on the one hand, we have
\begin{equation*}
\begin{split}
a_0 & = \frac{1}{2 ^{2 \alpha}}\left( \sum_{r=1}^{\beta}\frac{2}{3^{2\beta-r+2}}+\sum_{r=1}^{\beta}\sum_{s=0}^{r-1}\frac{4}{3^{2\beta-r-s+2}}  \right) \\
&=\frac{1}{2 ^{2 \alpha}}\left( \sum_{r=1}^{\beta}\frac{2}{3^{2\beta-r+2}}+\sum_{r=1}^{\beta}\frac{4}{3^{2\beta-r+2}}\cdot \frac{3^ r-1}{2}  \right) \\ 
&= \frac{1}{2 ^{2 \alpha}}\left( \sum_{r=1}^{\beta}\frac{2}{3^{2\beta-r+2}}+\sum_{r=1}^{\beta}\left(\frac{2}{3^{2\beta-2r+2}}-\frac{2}{3^{2\beta-r+2}}  \right) \right) \\
&= \frac{1}{2 ^{2 \alpha}}\cdot \frac{2}{3^{2\beta+2}}
\cdot \sum_{r=1}^{\beta}9^r=\frac{1}{2 ^{2 \alpha}}\cdot \frac{2}{3^{2\beta+2}} \cdot \frac{9^{\beta+1}-9}{8} 
\end{split}
\end{equation*}
Hence
\begin{equation}\label{eq: a0}
\begin{split}
a_0  = \frac{1}{2^{2\alpha+2}}\left(1-\frac{1}{3^{2\beta}} \right)\,.
\end{split}
\end{equation}
On the other hand, it is clear that the sum of the absolute values of the coefficients of $P$ is equal to the sum over $r$ and $s$ of the products 
of the sums of the absolute values of the coefficients of the polynomials 
$P_{1}$, $P_2$, $P_{3,r}$, and $P_{4,r,s}$. Note that the sum of the absolute values of the coefficients of $P_{1}$ or $P_2$ is $1$, whereas the sum of the absolute values of the coefficients of $P_{3,r}$ or $P_{4,r,s}$ is bounded by $2$. This yields the bound 
\begin{equation}\label{eq: bound}
\sum_{i=1}^{N} |a_i| \le \sum_{i=0}^{N} |a_i|\ \leq 1\cdot 1 \cdot \beta^ 2\cdot 2\cdot 2 = 4 \beta^2 \,.
\end{equation}

Putting~\eqref{eq: sumBk},  \eqref{eq: a0}, and~\eqref{eq: bound} toghether, it follows from Lemma~\ref{lem:Weil} that
\begin{equation}\label{eq: bound2}
\left| \sum_{k\in \F_\ell} \cB_{k,\ell}-\frac{\ell}{2^{2\alpha+2}}\left(1-\frac{1}{3^{2\beta}}\right) \right|\leq 4\beta^2\sqrt \ell\,.
\end{equation}
It is immediate that $\cB_{0,\ell}=a_0$. Observe that
\begin{equation*}
\begin{split}
\cB_{-1,\ell}&=\left( \frac{1}{2 ^{\alpha}}+ \frac{1}{2^\alpha}
\sum_{\chi\in \chi_{\ell,2^\alpha}^*}\chi(-1)\right)\frac{1}{2^\alpha} 
\\
& \qquad \qquad \cdot\sum_{r=1}^{\beta} \sum_{s=0}^{r-1} \frac{2+\vartheta_s}{3^{2\beta-s-r+2}}\\
& \qquad \qquad \qquad\cdot
\left(2+\vartheta_r + 3\sum_{\chi\in \chi_{\ell,3^{\beta-r}}^*}\chi(-1)-
\sum_{\chi\in \chi_{\ell,3^{\beta-r+1}}^*}\chi(-1)\right)\\
&=(2^\alpha-2)a_0\,.
\end{split}
\end{equation*}
For the last equality we have used that $\chi(-1)=1$ if the order of $\chi$ is a power of $3$, and the equality
$$
\sum_{\chi\in \chi_{\ell,2^\alpha}^*}\chi(-1)=2^\alpha-3\,.
$$
Combining~\eqref{eq:K and M} and~\eqref{eq: bound2}, we obtain
$$
\left|\# \Kl^*- \frac{\ell-2^\alpha+1}{2 ^{2 \alpha+2}}\left(1-\frac{1}{ 3 ^{2 \beta}}\right)\right| \leq 4 \beta ^2\sqrt \ell\,.
$$
Recalling~\eqref{eq:KK} we obtain 
$$
\left|\# \Kl- \frac{\ell}{2 ^{2 \alpha+2}}\left(1-\frac{1}{ 3 ^{2 \beta}}\right)\right|
 \leq 4 \beta ^2\sqrt \ell + 2 + \frac{2^\alpha-1}{2 ^{2 \alpha+2}}\left(1-\frac{1}{ 3 ^{2 \beta}}\right)\,.
$$
Since 
$$
2 + \frac{2^\alpha-1}{2 ^{2 \alpha+2}}\left(1-\frac{1}{ 3 ^{2 \beta}}\right)\leq \frac{33}{16}
$$
for $\alpha = 1,2, \ldots$, 
the result now follows. 

\section{Proof of Theorem~\ref{thm:Kl empty}}

We see from Corollary~\ref{cor:positive} that if $\#\Kl= 0$ then
$\ell  < 441 \cdot 2^{4\alpha } (\log x)^4$. So for $m$ in the 
representation $\ell-1= 2^\alpha 3^\beta m$ we have $m = O(\ell^{3/4}\log x) 
=  O(x^{3/4}\log x)$. 
Clearly for every $m =  O(x^{3/4}\log x)$ there are $O((\log x)^2)$ 
pairs on nonnegative integers $(\alpha, \beta)$ with $ 2^\alpha 3^\beta m+1 \le x$.
This concludes the proof. 

\section{Proof of Theorem~\ref{thm:vanish}}

Let us assume that $\beta\geq 1$ (otherwise the statement is immediate). Suppose that there exists $k\in \Kl$. Then Lemma~\ref{lemma: charac} implies that 
$$
\ordl k=3^ad\,,\qquad \ordl (-k^2-k)=3^be\,,
$$ 
where $1\leq a\leq \beta$, $0\leq b\leq a-1$, and $d,e$ are divisors of $m$. Note that $d$ must be nontrivial if $a=1$ and $e$ must be nontrivial if  $b=0$. Let $\Phi_m(X)$ denote the $m$th cyclotomic polynomial. Then $k$ is simultaneously a root of
$$
p_{a,d}(X):=\Phi_{3^ad}(X)% =x^{2\cdot 3^{a-1}}+x^{3^ {a-1}}+1 
\mand q_{b,e}(X):=\Phi_{3^bd}(-X^2-X)
$$
modulo $\ell$. This means that $\ell$ divides the resultant 
$$R_{a,b,d,e}:=\Res(p_{a,d},q_{b,e}),
$$

Note  that $p_{a,d}(X)$ and $q_{b,e}(X)$ have no roots in common. Indeed, let~$\zeta$ be a root of $p_{a,d}(X)$ (and so a root of unity of order $3^ad$) that is also a root of~$q_{b,e}(X)$. Then~$\zeta$ must satisfy that $-\zeta^ 2-\zeta=\eta $ or, equivalently, 
$$
\zeta+1=-\frac{\eta}{\zeta}\,,
$$
where $\eta$ is root of unity of order $3^ be$. Note that if a root of unity plus~1 is again a root of unity, then this root of unity is a primitive cubic root of unity. This is a contradiction with the fact that we can not have $a=1$ and $d=1$ simultaneously. Hence $R_{a,b,d,e}\ne 0$.

Furthermore, since all roots $\zeta$ of $p_{a,d}(X)$ have absolute value 1, we have
$$
\left|R_{a,b,d,e}\right| = \prod_{\zeta:~p_{a,d}(\zeta)=0} \left|q_{b,e}(\zeta) \right|
= \exp\(O(3^{a+b}\varphi(d)\varphi(e))\)\,,
$$
where the product runs over the $\zeta\in \C$ such that $p_{a,d}(\zeta) = 0$.
Note that if $\omega(t)$ is the number of distinct prime divisors of an integer~$t\ge 2$, then one has the inequality~$\omega(t)! \le t$.
Using Stirling's Formula, we derive $\omega(t) = O(\log t/\log (1 + \log t))$.  
Hence $R_{a,b,d,e}$ has at most 
$$
\omega(R_{a,b,d,e}) = O\(\frac{3^{a+b}\varphi(d)\varphi(e)}{a+b}\)
$$ 
distinct prime divisors. Hence
\begin{equation*}
\begin{split}
\# \cL_{\beta,m} & = O\(\sum_{a=1}^\beta\sum_{b =0}^{a-1}\sum_{d\mid m}\sum_{e\mid m}\frac{3^{a+b}\varphi(d)\varphi(e)}{a+b}\)\\
&= O\(\sum_{a=1}^\beta \sum_{b =0}^{a-1}\frac{3^{a+b}}{a}m^2\)= O\(\sum_{a=2}^\beta \frac{3^{2a}}{a}m^ 2\) = O\(3^{2\beta}m^ 2/\beta\). 
\end{split}
\end{equation*}

\section{Comments}

In addition to  $\#\cL_{0,m}=0$ of Theorem~\ref{thm:vanish} we also note that 
$\cL_{1,1} = \cL_{2,1} =\cL_{3,1} =\emptyset$.  Indeed, for $\cL_{1,1}$ the statement is immediate.
We now let $\cL_{a,b,d,e}$ denote the set of primes dividing $R_{a,b,d,e}$. One computes
$$
\cL_{2,1,1,1}=\{3\}, \qquad \cL_{3,2,1,1}=\{3,271\}, \qquad  \cL_{3,1,1,1}=\{3,271 \}\,. 
$$
It remains to note that $\cK_3=0$ and that $271$ is not of the 
form $2^\alpha 3^3+1$ for any $\alpha$. 

We now define $\ell_s$ as the smallest prime $\ell\equiv 1\pmod 3$
with $\cK_\ell = 0$ 
and $\omega(\ell-1)\ge s$ (if such prime exists), where, as before,
$\omega(t)$ denotes the number of distinct prime divisors of an integer~$t\ge 2$. 
From Theorem~\ref{thm:vanish} we expect that in fact 
 $\ell_s$ exists for any $s\ge 2$.  In Table~\ref{tab:ls} 
we
present some computational results which characterise the growth $\ell_s$.

\begin{table}[H]
  \centering
\begin{tabular}{|c|c|c|}
\hline
$s$  &  $\ell_s $  & Factorization of $\ell_s-1$\\
\hline
$3$  &  $31$         & $2\cdot 3\cdot 5$\\
$4$  &  $3121 $      & $2^4\cdot 3  \cdot 5\cdot 13$\\
$5$  &  $127681$    & $2^6\cdot 3\cdot 5 \cdot 7 \cdot 19$\\
$6$  &  $25858561$ & $2^9\cdot 3 \cdot5 \cdot7 \cdot13 \cdot37$\\
\hline
\end{tabular}
\caption{Values $\ell_s$ with $3\le s\le 6$}
\label{tab:ls}
\end{table}
We have not found $\ell_7$, but our computation shows that if  $\ell_7$ exists 
then $\ell_7 > 31\cdot 10^6$. 
On the other hand, combining the bound 
of  Theorem~\ref{thm:vanish} with the standard heuristic on the distribution
on primes,  one can derive a heuristic upper bound on $\ell_s$.

We remark that it is shown in~\cite{FGL} that if $k$ 
satisfies the conditions of Lemma~\ref{lemma: charac}, 
that is, $k \in \Kl$, then the rank $\rk\(D_{k,\ell}\)$ 
of the corresponding Demjanenko matrix satisfies
$$
\rk\(D_{k,\ell}\) = \frac{\ell-1}{2}\(1 - \frac{2}{M(k,\ell)}\), 
$$
where 
$$
M(k,\ell) := \lcm[\ordl(-k^2 -k), \ordl(k)].
$$
The resultant argument of~\cite{FGL} shows that 
$$
\min_{k \in \Kl} M(k,\ell) \to \infty
$$
as $\ell \to \infty$. This can easily be sharpened as
$$
\min_{k \in \Kl} M(k,\ell)  \ge c \sqrt{\log \ell}
$$
for an absolute constant $c>0$. In fact the same argument 
shows that for any real function $\psi(z)$ with $\psi(z)\to 0$ 
as $z\to \infty$,  all but $o(x/\log x)$ primes $\ell \le x$ we
have 
$$
\min_{k \in \Kl} M(k,\ell)  \ge \psi(\ell) \ell^{1/3}.
$$
Finally,  we remark that our approach allows to study the 
distribution of the values of $M(k,\ell) $ for every $\ell$.
%
%\section{Addendum?}

%For fixed integers $\beta \ge 0$ and $m\geq 1$ such that $(m,6)=1$, define $\cL_{\beta,m}$ to be the set of primes of the form $\ell=  2^\alpha 3^\beta m +1$, for some $\alpha > 0$, such that $\#\Kl>0$.  

%\begin{theorem}
%\label{thm:vanish} 
%For any fixed~$\beta\geq 0$ and $m\geq 1$ such that $(m,6)=1$, the set $\cL_{\beta,m}$ is finite. More precisely, if $\beta=0$, then $\#\cL_{\beta,m}=0$, and we have
%$$
%\# \cL_{\beta,m} =O\( 3^{2\beta}m^ 2/\beta \)
%$$
%for $\beta \ge 1$. 
%\end{theorem}

\section*{Acknowledgements}

The interest for the proportion of $k\in [2, \ell-2]$ for which $D_{k,\ell}$ is singular 
arose after a question of Kiran Kedlaya during a talk of the first author at the Workshop ``Frobenius distributions on curves'' held at CIRM, Luminy, in February 2014. The authors are grateful 
to CIRM for its support and hospitality. 
%The second author rapidly adverted a method to determine the 
%asymptotic proportion. Thanks to every one that made the conference possible.

During the preparation F.~Fit\'e was funded by the German Research Council via CRC 701, and partially supported by MECD project MTM2012-34611; I.~E.~Shparlinski was supported in part
by ARC grant DP130100237.


\begin{thebibliography}{McK-Sta}


\bibitem[FGL14]{FGL} 
F. Fit{\'e}, J. Gonz{\'a}lez, J.-C. Lario, 
`Frobenius distribution for quotients of Fermat curves of prime exponent', 
{\it Preprint\/}, 2014, 
(available from {\tt http://arxiv.org/abs/1403.0807}).

\bibitem[Gre80]{Gre80} R. Greenberg, `On the Jacobian variety of some algebraic curves', 
{\it Compos. Math.\/}, {\bf 42} (1980/81), 345--359.

\bibitem[IK04]{IwKow} H. Iwaniec and E. Kowalski,
{\it Analytic number theory\/}, Amer.  Math.  Soc.,
Providence, RI, 2004.

\bibitem[KR78]{KR78}
N. Koblitz and D. Rohrlich,
`Simple factors in the Jacobian of a Fermat curve',
{\it Canadian J. Math.\/}, {\bf 30} (1978), 1183--1205.

\end{thebibliography}
\end{document}